\documentclass[11pt,reqno]{amsart}
\usepackage[left=0.95in,right=0.95in, top=1.3in ,bottom=1.3in]{geometry}
\usepackage{bm, makecell}
\usepackage{float}
\usepackage{url}
\restylefloat{table}
\usepackage{amssymb,latexsym}
\usepackage{color}
\usepackage{booktabs}
\usepackage{graphicx}
\usepackage{epstopdf}
\usepackage{subfig}
\usepackage[mathcal]{euscript}
\usepackage{mathrsfs,amsmath}
\usepackage{enumerate}
\usepackage{hyperref}
\usepackage{amsfonts}
\usepackage{amssymb}
\usepackage{tikz}
\usetikzlibrary{shapes.geometric, arrows}
\usepackage{cite}
\usepackage[titletoc,title]{appendix}
\newtheorem{thm}{Theorem}[section]

\theoremstyle{definition}

\title[Dynamics of Zika infection with Gambusia Affinis]{Nonlinear Dynamic analysis of vector-host model for Zika infection with predatory fish Gambusia Affinis}
\author[S. Kumar]{Sachin Kumar}
\author[S. Jain]{Shikha Jain$^\ast$} \thanks{$^\ast$Corresponding Author}

\address{Sachin Kumar,
Assistant Professor,
Department of Mathematics,
University of Delhi,
New Delhi-110007,
India.}
\email{sachinambariya@gmail.com}

\address{Shikha Jain,
Research Scholar, 
Department of Mathematics,
University of Delhi,
New Delhi-110007,
India.}
\email{shikhajain01051990@gmail.com}



\numberwithin{equation}{section}
\allowdisplaybreaks

\begin{document}

\begin{abstract}
In the present paper, we study the dynamics of a nine compartmental vector-host model for Zika virus infection where the predatory fish Gambusia Affinis is introduced into the system to control the zika infection by preying on the vector. The system has six practically feasible equilibrium points where four of them are disease-free, and the rest are endemic. We discuss the existence and stability conditions for the equilibria. We find that when sexual transmission of zika comes to a halt then in absence of mosquitoes infection cannot persist. Hence, one needs to eradicate mosquitoes to eradicate infection. Moreover, we deduce that in the case of zika infection pushing the basic reproduction number below unity is next to impossible. Therefore, $\mathcal{O}_0$, the mosquito survival threshold parameter and $\mathcal{O}$, the mosquito survival threshold parameter with predation play a crucial role in getting rid of the infection in respective cases since mosquitoes cannot survive when these are less than unity. Sensitivity analysis shows the importance of reducing mosquito biting rate and mutual contact rates between vector and host. It exhibits the importance of increasing the natural mortality rate of vectors to reduce the basic reproduction number. Numerical simulation shows that when the basic reproduction number is close but greater than unity, the introduction of a small amount of predatory fish Gambusia Affinis can completely swipe off the infection. In case of high transmission or high basic reproduction number, this fish increases the susceptible human population and keeps the infection under control, hence, prohibiting it from becoming an epidemic.  

\end{abstract}

\maketitle

\section*{\textbf{Introduction}}
Zika virus was first isolated in 1947, from the blood of a sentinel rhesus macaque which was placed in the Zika Forest of Uganda for a study of yellow fever \cite{incubation}. Next, in January 1948, it was isolated from Aedes africanus mosquitoes, in the same forest. In 1952, the first human cases were detected in the United Republic of Tanzania and Uganda \cite{2overview, 3overview}. Few human cases were detected, mainly by serological methods, but such cases were scarce and the disease was considered to be benign. ZIKAV infection was regarded as largely asymptomatic or to cause only a mild illness with symptoms such as fever, rash, conjunctivitis, headache, and joint pain, etc.

  Zika virus is a Flavivirus, belonging to the family Flaviviridae. It is related to dengue, yellow fever, Japanese encephalitis, and West Nile Virus and is primarily transmitted by Aedes species of mosquitoes, particularly Aedes aegypti. These mosquitoes usually bite during the daytime. This virus is also passed on from a mother to its fetus during pregnancy, through breastfeeding, blood and plasma transfusion, and organ transplantation. Zika virus transmitted from a mother to its fetus during pregnancy may cause microcephaly and further congenital malformations in the infant, comprehensively known as congenital Zika syndrome. This virus is also spread through sexual contact though very few cases are reported so far. The diagnosis of Zika infection can only be established by laboratory tests of blood or other body fluids, such as semen or urine. In rare cases, zika may cause Guillain-Barre syndrome which is an uncommon nervous system sickness in which the own immune system of the victim starts damaging the nerve cells, causing muscle weakness, and occasionally, paralysis. Deaths due to zika infection are very rare. Since the primary transmitter of ZIKAV is the Aedes mosquito, people are advised to avoid mosquito bites by wearing full-coverage clothes treated with permethrin, using mosquito repellents, and using mosquito nets, etc. Using mosquito control strategies such as eliminating stagnant water, repairing septic tanks, and using insecticides to kill mosquitoes are preferred. Another way to control the mosquito population is by using genetically modified mosquitoes. The genetically modified male mosquitoes then mate with the pest female mosquitoes to produce the off-springs which die before reaching adulthood. But the experiment had uncertain results \cite{geneticex}. Moreover, genetic modification is an expensive and rigorous process with uncertain outcomes. Hence, the need of relying on an effective method for mosquito eradication is increased. 
  
  In the present paper, we try to analyze how the introduction of Gambusia Affinis into the system affects the infected human population. Gambusia Affinis is a freshwater fish commonly known as mosquitofish since they are larvivorous. Mosquitofish are viewed as the measure to control mosquito larvae as an appealing and safer substitute to the use of insecticides. In some areas, natural dispersal has caused the range extensions far from sites where they were initially introduced. Warm, shallow, and slow-moving waters are the optimal environment for G. affinis with dense vegetation, ample food organisms, and high mineral content. Precisely, they flourish in fresh and brackish water. They are mainly found in ponds, lakes, rivers, ditches, creeks, and springs.  They cannot tolerate temperatures below 4 degrees Celsius \cite{gambusiapredatoryuse}. A lot of prey-predator research papers have been published for using predatory fishes as a mosquito control measure \cite{gam1, gam2, gam3, gam4, gam5, gam6, gam7, gam8, gam9, gam10} but in the present study, we try to analyze the direct effect of introducing G. affinis to the disease model. A lot of research has been already done on zika virus dynamics \cite{zika1, zika2, zika3, zika4, zika5, zika7}. The papers on using this fish as an aid to mosquito control and the paper \cite{ppm}, where the prey-predator model is used for pest control are the motivation to this paper. In the present paper, we couple a vector-host disease model with a prey-predator model.
  
  The present paper consists of four sections: The section \ref{Section2} discusses the model formulation and basic properties like positivity and boundedness of solutions. Section \ref{sec2} gives the practically feasible equilibria of the system and explores their existence, stability conditions, and calculate the basic reproduction number. Section \ref{sec3} is dedicated to the numerical simulation and sensitivity analysis. Section \ref{sec4} summarizes all the results and conclusions.

\section{\textbf{Model formulation and basic properties}}\label{Section2}
The model subdivides the human population into four classes containing susceptible, incubating, infectious, and recovered individuals denoted by $S_h,\; E_h, \; I_h$ and $R_h$ respectively. Similarly, the vector population is subdivided into four classes containing immature stages (eggs, larvae, pupae), susceptible, incubating, and infectious individuals denoted by $m_q,\;S_v,\;E_v$ and $I_v$ respectively. The vector population does not have a recovery class, since their life span is short and the infective period ends with their death. The last class is of Gambusia Affinis, denoted by $G$. Moreover, $N_h = S_h + E_h + I_h + R_h$ and $N_v =  S_v + E_v + I_v$ are respectively, the total host and adult vector populations at time t. 

We assume the populations to be homogeneous, that is, all the individuals in a given class are identically infectious. We do not consider the vertical transmission in mosquito population because according to a research of University of Texas Medical Branch, only 1 per 290 offspring of an infected mosquito is infected which may be a male or a female mosquito but vertical transmission in humans is considered with $p$ and $q$ be the fraction of babies born infected which belong to the exposed and infected class respectively. We do not consider the sexual transmission of zika in humans as very few cases are reported so far. The susceptible human population is increased by the recruitment rate $\Lambda_h$ and $B_h$ is the human birth rate. We assume $d_h$ and $d_v$ to be the natural mortality rates, whereas $d_i$ and $d_{iv}$ are the disease-induced death rates for humans and adult female mosquitoes respectively while $d_q$ is the mortality rate of the immature stages. We assume $\epsilon$ is the oviposition rate, $K_q$ is the maximum carrying capacity of the immature stages, $\theta$ is the recovery rate from infected to recovered class, $k$ is the fraction of viable eggs that hatch into larvae, $\alpha$ is the transition rate from immature to mature stage and $f$ is the fraction of female mosquitoes, $\phi$ is predation rate of mosquitoes by the predator and $\psi$ is the additional growth rate of the predator due to predation. The man biting rate of mosquitoes is assumed to be $a_v$ and $\rho$ is the intrinsic growth rate of the predator, $C_{vh}$ and $C_{hv}$ are the transmission probabilities from an infectious mosquito to a susceptible human per bite and from an infected human to a susceptible mosquito per bite respectively, whereas $1/\nu_h$ and $1/\nu_v$ are the extrinsic incubation period in humans and mosquitoes respectively.

 All the above assumptions lead to the formation of following mathematical model:  
\begin{align}\label{eq:main}
\frac{dS_h}{dt} &= \Lambda_h - pB_hE_h-qB_hI_h-\frac{C_{vh}a_v I_v}{N_h}S_h - d_h S_h,\nonumber\\
\frac{dE_h}{dt} &= p B_h E_h + \frac{C_{vh}a_v I_v}{N_h}S_h -d_h E_h -\nu_h E_h,\nonumber\\
\frac{d I_h}{dt} &= qB_hI_h+\nu_h E_h-(\theta+d_h+d_i) I_h,\nonumber\\
\frac{d R_h}{dt} &=\theta I_h -  d_h R_h,\nonumber\\
\frac{d m_q}{d t} &= \epsilon k N_v \left(1-\frac{m_q}{K_q}\right) - (\alpha+d_q+\phi G)m_q,\nonumber\\
\frac{dS_v}{dt} &=f \alpha m_q - \frac{C_{hv}a_v I_h}{N_v}S_v-d_v S_v,\nonumber\\
\frac{dE_v}{dt} &= \frac{C_{hv}a_v I_h}{N_v}S_v -d_v E_v - \nu_v E_v,\nonumber\\
\frac{dI_v}{dt} &= \nu_v E_v -(d_v+d_{iv})I_v,\nonumber\\
\frac{dG}{dt} &= (\rho+\psi m_q)G\left(1-\frac{G}{K_x}\right).\nonumber\\
\end{align}

\subsection{\textbf{Positivity and boundedness of solutions}}

Since our system deals with host and vector populations which can not be negative, it is trivial to deduce that all the variables are always non-negative and the solutions of system \eqref{eq:main} always remain positive with positive initial conditions in the bounded region defined by
\begin{equation*}
\begin{split}
 \Delta =& \biggl\{(S_h, E_h, I_h, R_h, m_q, S_v, E_v, I_v, G)\in\mathbb{R}_+^{9}: N_h(t)\leqslant\frac{\Lambda_h}{d_h}, m_q(t) \leqslant K_q,\\ & N_v(t)\leqslant\frac{f \alpha K_q}{d_v}, 
  G(t)\leqslant K_x\biggr\}.
 \end{split}
\end{equation*}
and with the direct algebraic calculations we get the next result. 
\begin{thm}\label{th1}
For the non-negative initial conditions, the solutions of the system \eqref{eq:main} are positive when $t\geqslant0$ and the region $\Delta$ is positively invariant.
\end{thm}

\section{\textbf{Equilibria and stability analysis}}\label{sec2}
The system \eqref{eq:main} has following equilibrium points of the form

 $\mathcal{E}=\{S_h, E_h, I_h, R_h, m_q, S_v, E_v, I_v, G\}$:
\begin{itemize}
\item $\mathcal{E}_1= \{\Lambda_h/d_h , 0,0,0,0,0,0,0,0\}$.
\item $\mathcal{E}_2= \{\Lambda_h/d_h , 0,0,0,0,0,0,0,K_x\}$.
\item $\mathcal{E}_3= \{\Lambda_h/d_h , 0,0,0,K_q \left(1-\frac{1}{\mathcal{O}_0}\right), \frac{f \alpha K_q}{d_v}\left(1-\frac{1}{\mathcal{O}_0}\right),0,0,0\}$ where $\mathcal{O}_0 = \epsilon k f \alpha/{d_v(\alpha+d_q)}$ which is mosquito survival threshold without predation or the average number of female mosquitoes produced per adult female mosquito when predator is absent \cite{arias, chitnis}. It is equivalent to the {\it basic offspring number} derived in \cite{thpr}. Clearly, $\mathcal{E}_3$ biologically exists when $\mathcal{O}_0>1$.
\item $\mathcal{E}_4= \{\Lambda_h/d_h , 0,0,0,K_q \left(1-\frac{1}{\mathcal{O}}\right), \frac{f \alpha K_q}{d_v}\left(1-\frac{1}{\mathcal{O}}\right),0,0,K_x\},$ where $\mathcal{O}$ is the mosquito survival threshold in the presence of predator and equals $\epsilon k f \alpha/{d_v(\alpha+d_q+ \newline \phi K_x)}$. Again, $\mathcal{E}_4$ biologically exists when $\mathcal{O}>1$.
\item $\mathcal{E}_5= \{ S_h^5, E_h^5, I_h^5, R_h^5, 0, 0, 0, 0, 0\}$ where
\begin{align*}
\begin{split}
S_h^5=&\frac{(\theta - q B_h + d_h + d_i) (-p B_h + d_h +\nu_h)\Lambda_h}{p q B_h ^2 d_h + (\theta + d_h + d_i) (d_h + \lambda_h) (d_h + \nu_h) - B_h d_h (p \theta + (p + q)d_h + p d_i + q( \lambda_h +\nu_h))},\\
 E_h^5 =& \frac{(\theta - q B_h + d_h + d_i)\lambda_h \Lambda_h}{p q B_h ^2 d_h + (\theta + d_h + d_i) (d_h + \lambda_h) (d_h + \nu_h) - B_h d_h (p \theta + (p + q)d_h + p d_i + q( \lambda_h +\nu_h))},\\
I_h^5 =& \frac{\lambda_h \nu_h \Lambda_h}{p q B_h ^2 d_h + (\theta + d_h + d_i) (d_h + \lambda_h) (d_h + \nu_h) - B_h d_h (p \theta + (p + q)d_h + p d_i + q( \lambda_h +\nu_h))},\\
R_h^5=& \frac{\theta \lambda_h \nu_h \Lambda_h}{d_h(p q B_h ^2 d_h + (\theta + d_h + d_i) (d_h + \lambda_h) (d_h + \nu_h) - B_h d_h (p \theta + (p + q)d_h + p d_i + q( \lambda_h +\nu_h)))}.
 \end{split}
\end{align*}
where $\lambda_h= c_{vh} a_v I_v/N_h$ and $\lambda_v=c_{hv} a_v I_h/N_v$ are the force of infection on humans due to mosquitoes and force of infection on mosquitoes due to humans respectively.
\item $\mathcal{E}_6= \{S_h^6=S_h^5,E_h^6=E_h^5,I_h^6=I_h^5,R_h^6=R_h^5,0,0,0,0,K_x\}$.
\item $\mathcal{E}_7= \{S_h^7=S_h^5,E_h^7=E_h^5,I_h^7=I_h^5,R_h^7=R_h^5, m_q^7, S_v^7, E_v^7, I_v^7,0\}$ where
\begin{align*}
m_q^7 = & \frac{K_q ((f k \alpha \epsilon - (\alpha + d_q) d_v) (d_v + \lambda_v) (d_v + \nu_v) + d_{iv} B_1)}{f k \alpha \epsilon ((d_v + \lambda_v)d_v + \nu_v) + d_{iv} d_v +\lambda_v + \nu_v))},\\
S_v^7 = & \frac{K_q ((f k \alpha \epsilon - (\alpha + d_q) d_v) (d_v + \lambda_v) (d_v + \nu_v) + d_{iv} B_1)}{ k \epsilon (d_v + \lambda_v)((d_v + \lambda_v)d_v + \nu_v) + d_{iv} d_v +\lambda_v + \nu_v))},\\
E_v^7 = & K_q \lambda_v \left( \frac{f \alpha}{(d_v + \lambda_v) (d_v + \nu_v)} - \frac{(\alpha + d_q) (d_{iv} + d_v)}{k \epsilon ((d_v + \lambda_v) (d_v + \nu_v) + d_{iv}( d_v + \lambda_v+ \nu_v))}\right),\\
I_v^7=  & \frac{K_q \lambda_v \nu_v((f k \alpha \epsilon - (\alpha + d_q) d_v) (d_v + \lambda_v) (d_v + \nu_v) + d_{iv} B_1)}{ k \epsilon (d_{iv}+d_v) (d_v+\nu_v) (d_v + \lambda_v)((d_v + \lambda_v)d_v + \nu_v) + d_{iv} d_v +\lambda_v + \nu_v))},
\end{align*}
and 
\begin{align*}
B_1=((f k \alpha \epsilon - (\alpha + d_q) d_v)( d_v + 
\lambda_v) + (f k \alpha \epsilon - (\alpha + d_q)(d_v + \lambda_v)) \nu_v.
\end{align*}
\item $\mathcal{E}_8= \{S_h^8=S_h^5,E_h^8=E_h^5,I_h^8=I_h^5,R_h^8=R_h^5, m_q^8, S_v^8, E_v^8, I_v^8,K_x\}$ where 
\begin{align*}
m_q^8 = & \frac{K_q ((f k \alpha \epsilon - (\alpha + d_q +\phi K_x) d_v) (d_v + \lambda_v) (d_v + \nu_v) + d_{iv} B_2)}{f k \alpha \epsilon ((d_v + \lambda_v)d_v + \nu_v) + d_{iv} d_v +\lambda_v + \nu_v))},\\
S_v^8 = & \frac{K_q ((f k \alpha \epsilon - (\alpha + d_q + \phi K_x) d_v) (d_v + \lambda_v) (d_v + \nu_v) + d_{iv} B_2)}{ k \epsilon (d_v + \lambda_v)((d_v + \lambda_v)d_v + \nu_v) + d_{iv} d_v +\lambda_v + \nu_v))},\\
E_v^8 = & \frac{K_q \lambda_v((f k \alpha \epsilon - (\alpha + d_q + \phi K_x) d_v) (d_v + \lambda_v) (d_v + \nu_v) + d_{iv} B_2)}{ k \epsilon (d_v+\nu_v) (d_v + \lambda_v)((d_v + \lambda_v)d_v + \nu_v) + d_{iv} d_v +\lambda_v + \nu_v))},\\
I_v^8=  & \frac{K_q \lambda_v \nu_v((f k \alpha \epsilon - (\alpha + d_q + \phi K_x) d_v) (d_v + \lambda_v) (d_v + \nu_v) + d_{iv} B_2)}{ k \epsilon (d_{iv}+d_v) (d_v+\nu_v) (d_v + \lambda_v)((d_v + \lambda_v)d_v + \nu_v) + d_{iv} d_v +\lambda_v + \nu_v))},
\end{align*} and 
\begin{align*}
B_2=((f k \alpha \epsilon - (\alpha + d_q + \phi K_x) d_v) (d_v + 
\lambda_v) + (f k \alpha \epsilon - (\alpha + d_q +\phi K_x)(d_v + \lambda_v)) \nu_v.
\end{align*}
\end{itemize}
In the next section, we calculate the basic reproduction number before discussing the existence and stability of the above equilibria.  
\subsection{Reproduction number}
The basic reproduction number or herd immunity threshold is the average number of secondary infections caused per infective during the period of infection in a population consisting of susceptible entirely. We calculate the reproduction number $\mathfrak{R}_0$ for the system \eqref{eq:main} using next generation matrix method \cite{Dkm1, Dkm2} and the corresponding transmission matrix $T$ and transition matrix $\Sigma$ are 
\begin{align*}
T=&\begin{bmatrix}
0 & 0 & 0 & a_v C_{vh}\\
0 &0 & 0 & 0\\
0 & a_v C_{hv} & 0 & 0\\
0 & 0 & 0 & 0
\end{bmatrix},\\
\Sigma=& \begin{bmatrix}
 p B_h -d_h-\nu_h & 0 & 0 & 0\\
 \nu_h & q B_h-\theta-d_h-d_i & 0 & 0\\
 0 & 0 & -d_v -\nu_v & 0\\
 0 & 0 & \nu_v & d_v +d_{iv}
\end{bmatrix},
\end{align*}
and
\begin{align}\label{sigmain}
 -\Sigma^{-1} = \begin{bmatrix}
\frac{1}{d_h+\nu_h-p B_h} & 0 & 0 & 0\\
\frac{\nu_h }{(\theta-q B_h +d_h+d_i)(d_h+\nu_h-p B_h)} &\frac{ 1}{(\theta-q B_h +d_h+d_i)}& 0 & 0\\
0 & 0 & \frac{1}{(d_v+\nu_v)} & 0\\
0 & 0 & \frac{\nu_v}{(d_{iv}+d_v)(d_v+\nu_v)} & \frac{1}{(d_{iv}+d_v)} 
\end{bmatrix}.
\end{align}
Reproduction number is the dominant eigenvalue of the matrix $-T\Sigma^{-1}$, hence calculations result as follows: 
\begin{equation}\label{rp}
\mathfrak{R}_0 = \sqrt{\frac{a_v^2 C_{vh} C_{hv}\nu_h \nu_v}{(q B_h-\theta-d_h-d_i)( p B_h -d_h-\nu_h)(d_{iv}+d_v)(d_v+\nu_v)}}.
\end{equation}
Moreover, simple calculations show that 
\begin{align*}
\mathfrak{R}_0 = \sqrt{\mathfrak{R}_0^h\, \mathfrak{R}_0^v}
\end{align*} where $\mathfrak{R}_0^h$ and $ \mathfrak{R}_0^v$ are the basic reproduction number for host and vector population respectively, given by the following
\begin{align*}
\mathfrak{R}_0^h =& \frac{a_v C_{vh}\nu_h}{(q B_h-\theta-d_h-d_i)( p B_h -d_h-\nu_h)},\\
\mathfrak{R}_0^v =& \frac{a_v C_{hv}\nu_v}{(d_{iv}+d_v)(d_v+\nu_v)}
\end{align*}

The element $(-\Sigma^{-1})_{ij}$ represents the expected time spent by an individual in the state i  during the epidemiological course, who presently has state j \cite{Dkmbook}. Since time is non-negative, from \eqref{sigmain} $(\theta-q B_h +d_h+d_i)$ and $(d_h+\nu_h-p B_h)$ are positive.
\subsection{Existence and stability of equilibria}
The first four equilibrium points are disease free equilibrium points and correspond to the force of infection $\lambda_h=0$ and $\lambda_v=0$. The equilibrium points $\mathcal{E}_3$ and $\mathcal{E}_4$ exist only when their corresponding mosquito survival thresholds are greater than unity. Moreover, it can be clearly seen that $\mathcal{O}_0> \mathcal{O}$, that is, predation decreases the mosquito threshold parameter. The equilibrium points $\mathcal{E}_5$ and $\mathcal{E}_6$ correspond to the force of infection $\lambda_v=0$ and 
\begin{align*}
\lambda_h=- \frac{d_h (\theta-q B_h+d_h+d_i)(d_h+\nu_h- p B_h)}{d_h (\theta-q B_h+d_h+d_i)+(\theta+d_h)\nu_h}
\end{align*} 
which is negative. Since force of infection can never be negative, the equilibrium points $\mathcal{E}_5$ and $\mathcal{E}_6$ are practically impossible to exist. This means that infection in human population can not persist without infection in mosquito population. The equilibrium points $\mathcal{E}_7$ and $\mathcal{E}_8$ correspond to the force of infection $\lambda_h>0$, $\lambda_v>0$ and $\lambda_h>0$, $\lambda'_v>0$, where  these are the force of infections in absence and presence of predator Gambusia Affinis respectively. Moreover, $\mathcal{E}_7$ and $\mathcal{E}_8$ exist when $p q B_h ^2 d_h + (\theta + d_h + d_i) (d_h + \lambda_h) (d_h + \nu_h)/ B_h d_h (p \theta + (p + q)d_h + p d_i + q( \lambda_h +\nu_h))>1$ as well as mosquito survival threshold $\mathcal{O}_0$ and mosquito survival threshold with predation $\mathcal{O}$ are greater than unity respectively.

Now, we discuss the stability of the equilibrium points. The Jacobian matrix of \eqref{eq:main} is as follows
\begin{align}\label{Jacobian}
J = \begin{bmatrix}
\frac{-C_{v h}I_v}{N_h}-d_h & -p B_h & -q B_h & 0 & 0 & 0 & 0 & -\frac{C_{vh}a_v S_h}{N_h} & 0\\
\frac{C_{v h}I_v}{N_h} & -B & 0 & 0 &0 &0 & 0 & \frac{C_{vh}I_v}{N_h} & 0\\
0 & \nu_h & -A & 0 & 0 & 0 & 0 & 0 & 0 \\
0 & 0 & \theta & -d_h & 0 & 0 & 0 & 0 & 0\\
0 & 0 & 0 & 0 & J_{55} & J_{5j} & J_{5j} & J_{5j} & -\phi m_q\\
0 & 0 & 0 & 0 & f \alpha & -d_v & 0 & 0 & 0\\
0 & 0 & a_v C_{hv} & 0 & 0 & 0 & -d_v - \nu_v & 0 & 0\\
0 & 0 &0 & 0 & 0 & 0 & \nu_v & -d_v-d_{iv} & 0\\
0 & 0 & 0 & 0 & \psi G\biggl(1-\frac{G}{K_x}\biggr) & 0 & 0 & 0 & J_{99}
\end{bmatrix}
\end{align}
where
\begin{align*}
J_{55}&=\frac{-\epsilon k N_v}{K_q}-(\alpha+d_q+\phi G),\\
J_{5 j}&=k \epsilon\biggl(1-\frac{m_q}{K_q}\biggr),\\
J_{99}&=\rho+\psi m_q \biggl(1-\frac{G}{K_x}\biggr)-(\rho+\psi m_q)\frac{G}{K_x}.\\
\end{align*}

The Jacobian matrix evaluated at the equilibrium point $\mathcal{E}_1$is given by
\begin{align}\label{JE_1}
J|_{\mathcal{E}_1} = \begin{bmatrix}
-d_h & -p B_h & -q B_h & 0 & 0 & 0 & 0 & -C_{vh} & 0\\
0 & -B & 0 & 0 &0 &0 & 0 & C_{vh} & 0\\
0 & \nu_h & -A & 0 & 0 & 0 & 0 & 0 & 0 \\
0 & 0 & \theta & -d_h & 0 & 0 & 0 & 0 & 0\\
0 & 0 & 0 & 0 & -\alpha-d_q & k \epsilon & k \epsilon & k \epsilon & 0\\
0 & 0 & 0 & 0 & f \alpha & -d_v & 0 & 0 & 0\\
0 & 0 & a_v C_{hv} & 0 & 0 & 0 & -d_v - \nu_v & 0 & 0\\
0 & 0 &0 & 0 & 0 & 0 & \nu_v & -d_v-d_{iv} & 0\\
0 & 0 & 0 & 0 & 0 & 0 & 0 & 0 & \rho 
\end{bmatrix},
\end{align}
where 
\begin{align*}
A =\, & \theta- q B_h+d_h+d_i,\\
B =\, & d_h+\nu_h -p B_h.
\end{align*}
The characteristic polynomial for \eqref{JE_1} is given by 
\begin{align}\label{CH:E_1}
\begin{split}
&(-x + \rho) (-x - d_h)^2 (x^2 + x \alpha - f k \alpha \epsilon + x d_q + x d_v + \alpha d_v + d_q d_v)((-A - x) (-B - x)\\& (-x - d_v - d_{i v}) (-x - d_v - \nu_v) - a_v^2 C_{vh} C_{hv}).
\end{split}
\end{align}
Clearly, from \eqref{CH:E_1}, $J|_{\mathcal{E}_1}$ has an eigenvalue with positive real part, hence $\mathcal{E}_1$ is not locally asymptotically stable or unstable independent of the reproduction number. Hence, we get the following result:
\begin{thm}
The equilibrium point $\mathcal{E}_1$ is not locally asymptotically stable irrespective of the basic reproduction number.
\end{thm}

The Jacobian matrix evaluated at the equilibrium point $\mathcal{E}_2$ is given by

\begin{align}\label{JE_2}
J|_{\mathcal{E}_2} = \begin{bmatrix}
-d_h & -p B_h & -q B_h & 0 & 0 & 0 & 0 & -C_{vh} & 0\\
0 & -B & 0 & 0 &0 &0 & 0 & C_{vh} & 0\\
0 & \nu_h & -A & 0 & 0 & 0 & 0 & 0 & 0 \\
0 & 0 & \theta & -d_h & 0 & 0 & 0 & 0 & 0\\
0 & 0 & 0 & 0 & -\alpha-d_q-\phi K_x & k \epsilon & k \epsilon & k \epsilon & 0\\
0 & 0 & 0 & 0 & f \alpha & -d_v & 0 & 0 & 0\\
0 & 0 & a_v C_{hv} & 0 & 0 & 0 & -d_v - \nu_v & 0 & 0\\
0 & 0 &0 & 0 & 0 & 0 & \nu_v & -d_v-d_{iv} & 0\\
0 & 0 & 0 & 0 & 0 & 0 & 0 & 0 & -\rho 
\end{bmatrix}
\end{align}
and the characteristic polynomial is given by
\begin{align}\label{CH:E2}
\begin{split}
&(-x - \rho) (-x - d_h)^2 (x^2 + x( \alpha + d_q +d_v +\phi K_x) - f k \alpha \epsilon + \alpha d_v + d_q d_v + \phi d_v K_x)((-A - x) \\&(-B - x) (-x - d_v - d_{i v}) (-x - d_v - \nu_v) - a_v^2 C_{vh} C_{hv}).
\end{split}
\end{align}
We use Ruth-Hurwitz criterion for determining the stability conditions for the last two factors of the equation \eqref{CH:E2}. Comparing the quadratic factor of \eqref{CH:E2} with the standard form $ax^2 + bx+c$, we get
\begin{align*}
a=&1,\\
b=& ( \alpha + d_q +d_v +\phi K_x),\\
c=&  - f k \alpha \epsilon + \alpha d_v + d_q d_v + \phi d_v K_x.
\end{align*} 
Clearly, $a$, $b$ are positive and $c$ is positive when $\frac{ \epsilon k f \alpha}{d_v(\alpha+d_q+\phi K_x)}=\mathcal{O}<1$. Comparing the next bi-quadratic polynomial with the standard form, that is, $a_0x^4+a_1x^3+a_2x^2+a_3x+a_4$, we get
\begin{align*}
a_0=& 1,\\
a_1=& A +B +2 d_v +d_{i v} +\nu_v,\\
a_2=& AB+2Ad_v+2Bd_v+d_v^2+A d_{i v}+ B d_{i v} +d_v d_{i v} +(A + B+ d_v +d_{i v})\nu_v,\\
a_3=& 2ABd_v+ (A+B)(d_v^2+d_v d_{i v}+ d_v \nu_v+d_{i v}\nu_v)+AB(d_{i v}+\nu_v),\\
a_4=& AB(d_v+d_{i v})(d_v + \nu_v)-a_v^2C_{vh}C_{hv}\nu_h \nu_v.
\end{align*}
Now, $a_0, \, a_1,\,a_2,\,a_3$ are positive and $a_4$ is positive when $\mathfrak{R}_0^2<1$. Further calculations result in the following Ruth-Hurwitz constants:
\begin{align}\label{c_1}
b_1&=\frac{AB(A+B)+(A+B+d_v+d_{i v})(A+B+d_v+\nu_v)(2d_v+d_{i v}+\nu_v)}{a_1},\\
b_2&= a_4,\\
c_1&= \frac{b_1a_3-a_1a_4}{b_1}.
\end{align}
$\therefore \mathcal{E}_2$ is a locally asymptotically stable critical point or a nodal attractor when $\mathcal{R}_0^2<1,\;\mathcal{O}<1$ and $c_1>0$. We summarize this result in the next theorem
\begin{thm}
The equilibrium point $\mathcal{E}_2$ is a nodal attractor when $\mathcal{R}_0^2<1,\;\mathcal{O}<1$ and $c_1>0$, where $c_1$ is the Ruth-Hurwitz constant given by \eqref{c_1}.
\end{thm}

Now, we discuss the stability of equilibrium point $\mathcal{E}_3$. 
The Jacobian matrix evaluated at $\mathcal{E}_3$ is given by
\begin{align}\label{JE_3}
J|_{\mathcal{E}_3}=\begin{bmatrix}
-d_h & -p B_h & -q B_h & 0 & 0 & 0 & 0 & -C_{vh} & 0\\
0 & -B & 0 & 0 &0 &0 & 0 & C_{vh} & 0\\
0 & \nu_h & -A & 0 & 0 & 0 & 0 & 0 & 0 \\
0 & 0 & \theta & -d_h & 0 & 0 & 0 & 0 & 0\\
0 & 0 & 0 & 0 & D & \frac{k \epsilon}{\mathcal{O}_0} & \frac{k \epsilon}{\mathcal{O}_0} & \frac{k \epsilon}{\mathcal{O}_0} & -\phi K_q (1-\frac{1}{\mathcal{O}_0})\\
0 & 0 & -a_v C_{h v} & 0 & f \alpha & -d_v & 0 & 0 & 0\\
0 & 0 & a_v C_{hv} & 0 & 0 & 0 & -d_v - \nu_v & 0 & 0\\
0 & 0 &0 & 0 & 0 & 0 & \nu_v & -d_v-d_{iv} & 0\\
0 & 0 & 0 & 0 & 0 & 0 & 0 & 0 & \rho +\biggl(1-\frac{1}{\mathcal{O}_0}\biggr)\psi K_q
\end{bmatrix},
\end{align}
where 
\begin{align*}
D=-\alpha-d_q-\frac{\alpha f k \epsilon(1-\frac{1}{\mathcal{O}_0})}{d_v}.
\end{align*}
The characteristic polynomial of \eqref{JE_3} is as follows
\begin{align}\label{CH:E_3}
\begin{split}
&\frac{1}{\mathcal{O}_0 d_v}(x+d_h)^2 \biggl(x-\rho-\biggl(1-\frac{1}{\mathcal{O}_0}\biggr)\biggr)((-A - x)(-B - x) (-x - d_v - d_{i v}) (-x - d_v - \nu_v)\\& - a_v^2 C_{vh} C_{hv}) (-f k \alpha \epsilon d_v +(x+d_v)(-f k \alpha \epsilon + f k \mathcal{O}_0 \alpha \epsilon + \mathcal{O}_0 x d_v + \mathcal{O}_0 \alpha d_v+ \mathcal{O}_0 d_q d_v  + \mathcal{O}_0 \phi d_v K_x)).
 \end{split}
\end{align} 
Clearly, \eqref{CH:E_3} has an eigenvalue with a positive real part. Hence, $\mathcal{E}_3$ is not locally asymptotically stable and the next result follows.
\begin{thm}
The equilibrium point $\mathcal{E}_3$ is never locally asymptotically stable.
\end{thm}

For the stability of $\mathcal{E}_4$, the Jacobian matrix evaluated at it, $J|_{\mathcal{E}_4}$ is given by
\begin{align}\label{JE_4}
\begin{bmatrix}
-d_h & -p B_h & -q B_h & 0 & 0 & 0 & 0 & -C_{vh} & 0\\
0 & -B & 0 & 0 &0 &0 & 0 & C_{vh} & 0\\
0 & \nu_h & -A & 0 & 0 & 0 & 0 & 0 & 0 \\
0 & 0 & \theta & -d_h & 0 & 0 & 0 & 0 & 0\\
0 & 0 & 0 & 0 & D & \frac{k \epsilon}{\mathcal{O}} & \frac{k \epsilon}{\mathcal{O}} & \frac{k \epsilon}{\mathcal{O}} & -\phi K_q (1-\frac{1}{\mathcal{O}})\\
0 & 0 & -a_v C_{h v} & 0 & f \alpha & -d_v & 0 & 0 & 0\\
0 & 0 & a_v C_{hv} & 0 & 0 & 0 & -d_v - \nu_v & 0 & 0\\
0 & 0 &0 & 0 & 0 & 0 & \nu_v & -d_v-d_{iv} & 0\\
0 & 0 & 0 & 0 & 0 & 0 & 0 & 0 & -\rho -\biggl(1-\frac{1}{\mathcal{O}}\biggr)\psi K_q
\end{bmatrix}
\end{align}
where
\begin{align*}
D = -\alpha-d_q-\frac{\alpha f k \epsilon(1-\frac{1}{\mathcal{O}})}{d_v}-\phi K_x.
\end{align*}
The characteristic polynomial for \eqref{JE_4} is given by
\begin{align}\label{CH:E_4}
\begin{split}
&\frac{1}{\mathcal{O}d_v}(x+d_h)^2 \biggl(x+\rho+\biggl(1-\frac{1}{\mathcal{O}}\biggr)\biggr)((-A - x)(-B - x) (-x - d_v - d_{i v}) (-x - d_v - \nu_v) - a_v^2 C_{vh} C_{hv})\\& (-f k \alpha \epsilon d_v +(x+d_v)(-f k \alpha \epsilon + f k \mathcal{O} \alpha \epsilon + \mathcal{O} x d_v + \mathcal{O} \alpha d_v+ \mathcal{O} d_q d_v  + \mathcal{O} \phi d_v K_x)).
 \end{split}
\end{align}
The first two factors of the polynomial in \eqref{CH:E_4} result in negative real roots as well as the third factor when $\mathcal{O}>1$ and $c_1>0$, as we proved for \eqref{CH:E2}. For the last quadratic factor, we again use Ruth-Hurwitz criteria for deducing the sign of real parts of the roots. We find that the all the non constant coefficients of the quadratic factor are positive and constant term is positive when $\mathcal{O}>2fk\alpha\epsilon /(fk\alpha\epsilon+(\alpha+d_q+\phi K_x)d_v)$. Hence, the equilibrium point $\mathcal{E}_4$ is locally asymptotically stable when these conditions are satisfied. Therefore, we get the following result:
\begin{thm}
The equilibrium point $\mathcal{E}_4$ is locally asymptotically stable when the following conditions are satisfied:
\begin{align*}
\mathcal{R}_0^2&<1,\\
c_1&>0,\\
\mathcal{O}&>\frac{2fk\alpha\epsilon}{fk\alpha\epsilon+(\alpha+d_q+\phi K_x)d_v}.
\end{align*}
\end{thm}
Next, we perform the numerical simulation to further get the precise conclusions.
\section{Numerical Simulation} \label{sec3} 
In the present section, we perform the numerical simulation using MATLAB. Table \ref{table:1} gives the numerical values of the parameters. 

\begin{table}[ht]
\caption{Parameters and their values}
\centering
\begin{tabular}{l c l c}
\hline \hline
Parameter & Symbol & Estimate & Source \\ [0.5ex]
\hline 
Human recruitment rate & $\Lambda_h$ & 19 & assumed\\
Human mortality rate & $d_h$ & 0.000039 & calculated\\
Human disease induced mortality rate & $d_i$ & 0.001 & \cite{di}\\
Human birth rate & $B_h$ & 0.019 & \cite{wb}\\
Vector mortality rate & $d_v$ & 0.042 & calculated\\
Vector disease induced mortality rate & $d_{iv}$ & 0.00001 & \cite{di}\\
Fraction of babies born in exposed class & $p$ & 0.07 & assumed\\
Fraction of babies born in infected class & $q$ & 0.07 & assumed\\
Man biting rate of mosquitoes & $a_v$ & 3.0425 & \cite{biting}\\
Extrinsic incubation period in humans & $1/\nu_h$ & 0.1429 & calculated\\
Extrinsic incubation period in mosquitoes & $1/ \nu_v$ & 0.25 & calculated\\
Recovery rate from infected to recovered class & $\theta$ & 0.02 & calculated\\
Oviposition rate & $\epsilon$ & 16.25 & calculated\\
Fraction of eggs hatched in to larvae & $k$ & 0.9 & assumed\\
Carrying capacity of the immature stages & $K_q$ & 5000 & assumed\\
Transition rate from immature to mature stages & $\alpha$ & 0.083 & calculated\\
Mortality rate of immature stages & $d_q$ & 0.048 & calculated\\
Predation rate & $\phi$ & 0.35 & calculated\\
Fraction of female mosquitoes & $f$ & 0.5 & assumption\\
Intrinsic growth rate of predator & $\rho$ & 0.37 & calculated\\
Predator's additional growth rate due to predation & $\psi$ & 0.05 & assumed\\
Carrying capacity of the predator & $K_x$ & 1000 & assumed\\ 
\hline
\end{tabular}
\label{table:1}
\end{table}
To calculate natural mortality rates, human life expectancy is considered to be 71.4 years \cite{le} and mosquito life expectancy is assumed to be 24 days as their life span can range from 2 to 4 weeks \cite{aedesls}. Zika incubation period in humans and mosquitoes ranges from 3-14 days and 2-7 days respectively \cite{incubation} and we choose the values to be 7 days and 4 days respectively. Moreover, we assume exactly half of the mosquito population are females. Zika virus remains detectable in the human bloodstream from 3-54 days \cite{recovery}, say for 50 days and hence, the recovery rate $\theta$ is 0.02. Oviposition rate is the total number of eggs produced per day by a female in her complete life span. Aedes mosquito can produce an average of 100-200 eggs per batch and a female can produce up to 5 batches of eggs during a lifetime \cite{aedesls}. We assume on average a female Aedes aegypti mosquito produces 130 eggs per batch and three batches in her life span, which results in an oviposition rate of 16.25 per day. An Aedes mosquito takes 8-15 days for becoming an adult from an egg \cite{aedesls}, say 12 days, and hence natural mortality rate for immature stages is 0.083. The maximum consumption in a day by one mosquito-fish has been observed to be 42-167 percent of its own body weight and males weight between 0.13-0.2 g while females weight 0.2-1.0 g \cite{Pyke}. We assume consumption to be 70 percent and average weight to be 0.5 g. Hence, the predation rate is 0.35. The growth rate of Gambusia Affinis is 20 percent of its dry weight \cite{Pyke} which results in the intrinsic growth rate of 0.37. Moreover, we are assuming a small man-made pond as the oviposition site and the place where we introduce Gambusia Affinis which can with-hold a maximum of 1000 fishes and we assume that it can also contain a maximum of 5000 immature stages of mosquitoes. Table \ref{table:1} enlists the parameters and the values used for the simulation.   
\subsection{Sensitivity Analysis}
For determining the best method to reduce infection in the human population, it becomes necessary to know the relative importance of different factors responsible for the prevalence and transmission of disease. The reproduction number, $\mathfrak{R}_0$ has a high impact on the initial disease transmission. Moreover, in previous sections, we analyzed that the mosquito threshold parameters $\mathcal{O}_0$ and $\mathcal{O}$ are of great importance in determining the existence of mosquito population and equilibrium points $\mathcal{E}_7$ and $\mathcal{E}_8$. Hence, we use sensitivity analysis to discover the parameters having a high impact on reproduction number and mosquito threshold parameters. The normalized forward sensitivity index of a variable, x, that depends differentiably on a parameter, p, is deﬁned as:
\begin{align*}
\Gamma^x_p = \frac{\partial x}{\partial p} \frac{p}{x}
\end{align*}
We evaluate the sensitivity indices at the baseline parameter values given in Table \ref{table:1} and the resulting sensitivity indices of the $\mathfrak{R}_0^2$ to the different parameters of the model are shown in table \ref{table:2}

\begin{table}[ht]
\caption{Sensitivity indices of $\mathfrak{R}_0^2$ to parameters of the model, evaluated at the baseline parameter values given in Table \ref{table:1}.}
\centering
\begin{tabular}{c c | c c}
\hline \hline
Parameter & Sensitivity Index & Parameter & Sensitivity Index\\
\hline \hline
$a_v$ & 2 & $d_{iv}$ & -0.000238039\\
$d_v$ & -1.1436 & $K_q$ & 0\\
$\theta$ & -1.01476 & $\alpha$ & 0\\
$c_{vh}$ & 1  & $\epsilon$ & 0\\
$c_{hv}$ & 1 & k & 0\\
$\nu_v$ & 0.123836 & $d_q$ & 0\\
$B_h$ & 0.0768739 & $\phi$ & 0\\
$q$ & 0.0674819 & $f$ & 0\\
$d_i$ & -0.0507382 & $\rho$ & 0\\
$p$ & 0.00939206 & $\psi$ & 0\\
$\nu_h$ & -0.00911665 & $K_x$ & 0\\
$d_h$ & -0.0022542 & $\Lambda_h$ & 0\\
\hline
\end{tabular}
\label{table:2}
\end{table}

\begin{table}[ht]
  \caption{Sensitivity indices of the mosquito threshold parameters $\mathcal{O}_0$ and $\mathcal{O}$, evaluated at the parameter values listed in Table \ref{table:1}.}
\makebox[\textwidth][c]{
    \begin{tabular}{cc|cc}
\toprule
  \multicolumn{2}{c|}{Analysis of $\mathcal{O}_0$} & \multicolumn{2}{c}{Analysis of $\mathcal{O}$} \\
\midrule
Parameter & Sensitivity Index & Parameter & Sensitivity Index\\
 \\ \hline
$\epsilon$ & 1 & $\epsilon$ & 1\\
$k$ & 1 & $k$ & 1\\
$f$ & 1 & f & 1\\
$d_v$ & -1 & $d_v$ & -1\\
$\alpha$ & 0.366412 & $\alpha$ & 0.999763\\
$d_q$ & -0.366412 & $d_q$ & -0.000137092\\
 & & $\phi$ & -0.999626\\
  & & $K_x$ & -0.999626\\
\bottomrule
\end{tabular}%
}
\label{table:3}%
\end{table}%

Table \ref{table:2} enlists the sensitivity indices of $\mathfrak{R}_0^2$ and hence shows the most important parameters for reducing infection in the human population. Most parameters are in accordance with the intuition like mosquito biting rate is directly related to the spread of infection. What surprises here is the higher sensitivity indices of mosquito mortality rate and human recovery rate even more than the mutual contact rates between host and vector. This suggests that if we could do something which can increase these two parameters would be very beneficial in controlling the infection. But till now there is no treatment for zika, hence one can not do anything with the recovery rate. In this model, we introduce Gambusia Affinis in the vector host system as a predator which preys on the vector and hence helps in decreasing their population size. 

Table \ref{table:3} enlists the parameters with non-zero sensitivity indices and shows the effect of concerned parameters on mosquito threshold parameters. Since these threshold parameters play a crucial role in determining the existence of mosquito population and hence the endemic equilibrium points. Clearly, $\mathcal{O}<\mathcal{O}_0$, therefore reducing the former to less than unity is much easier. The sensitivity indices of most of the parameters agree with an intuitive expectation. The outcome that perhaps need some clarification is that for the transition rate from immature to mature stages, $\alpha$. In the presence of G. Affinis, the dependency on $\alpha$ increases instead of decreasing and this may be explained as the higher transition rate would make it difficult for G. Affinis to devour the immature stages before the transition and hence adding to the mosquito survival threshold with predation. The addition of this predatory fish also lessens the dependability of mosquito threshold parameter on the natural mortality rate of immature stages since it leads to the forced mortality of immature stages by preying on them. The predation rate, $\phi$, and Carrying capacity of the predator, $K_x$ also have sensitivity indices approaching unity in a negative sense. Both of these parameters can be easily controlled and hence can be of great help in reducing mosquito threshold parameter with predation.

\subsection{Simulation} 
For simulation, we use the parameter base values from Table \ref{table:1} and first we try to calculate our reproduction number.
On substituting the values from Table \ref{table:1} in \eqref{rp}, we get
\begin{align*}
\mathfrak{R}_0= 98.2814 \sqrt{c_{vh} c_{hv}}.
\end{align*}

\begin{figure}[ht!]
\begin{center}
\subfloat[$G_0=0, c_{vh}=0.2, c_{hv}=0.25, a_v=0.25$]{\includegraphics[scale=0.5]{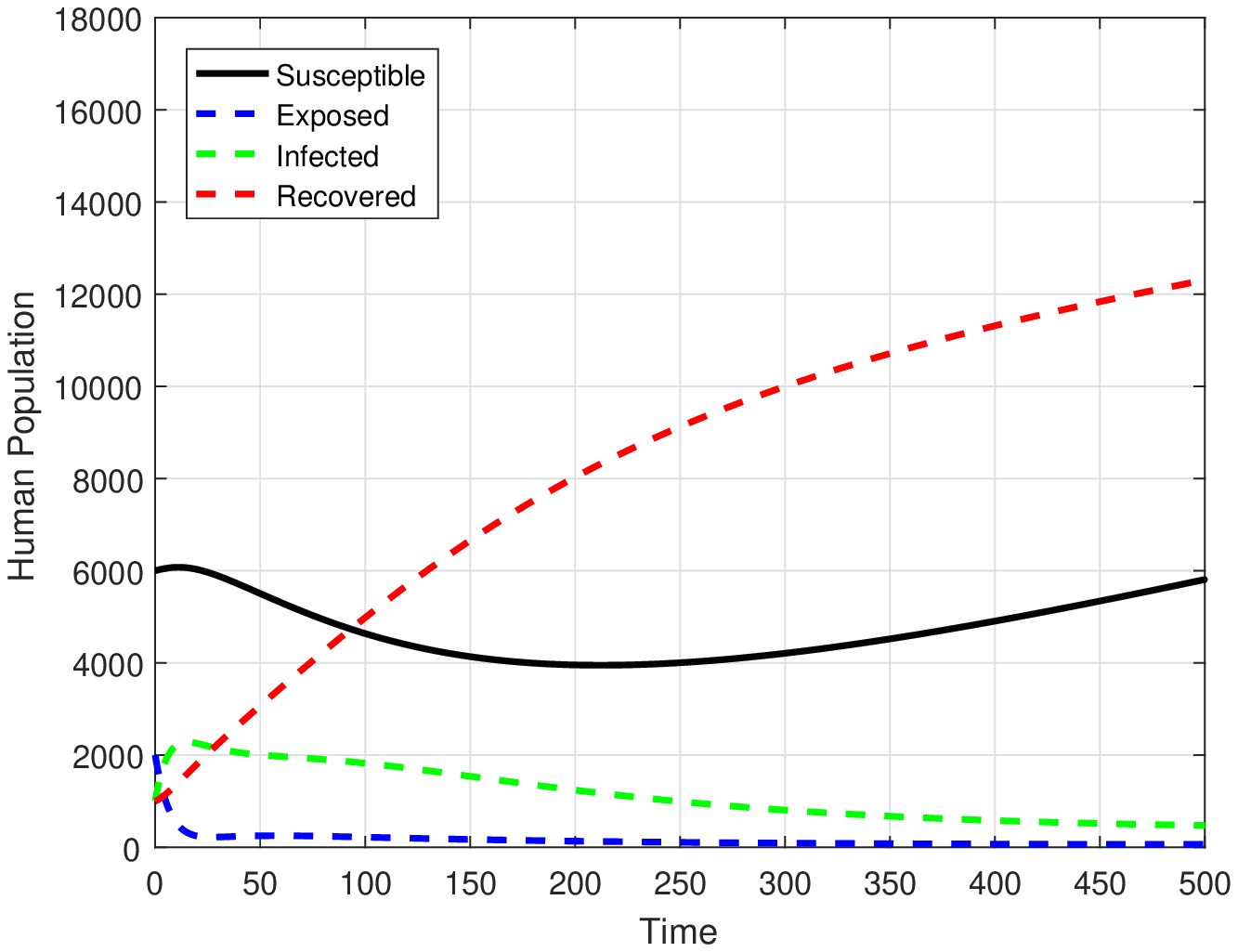} \label{a}}
\subfloat[$G_0=20, c_{vh}=0.2, c_{hv}=0.25, a_v=0.25$]{\includegraphics[scale=0.5]{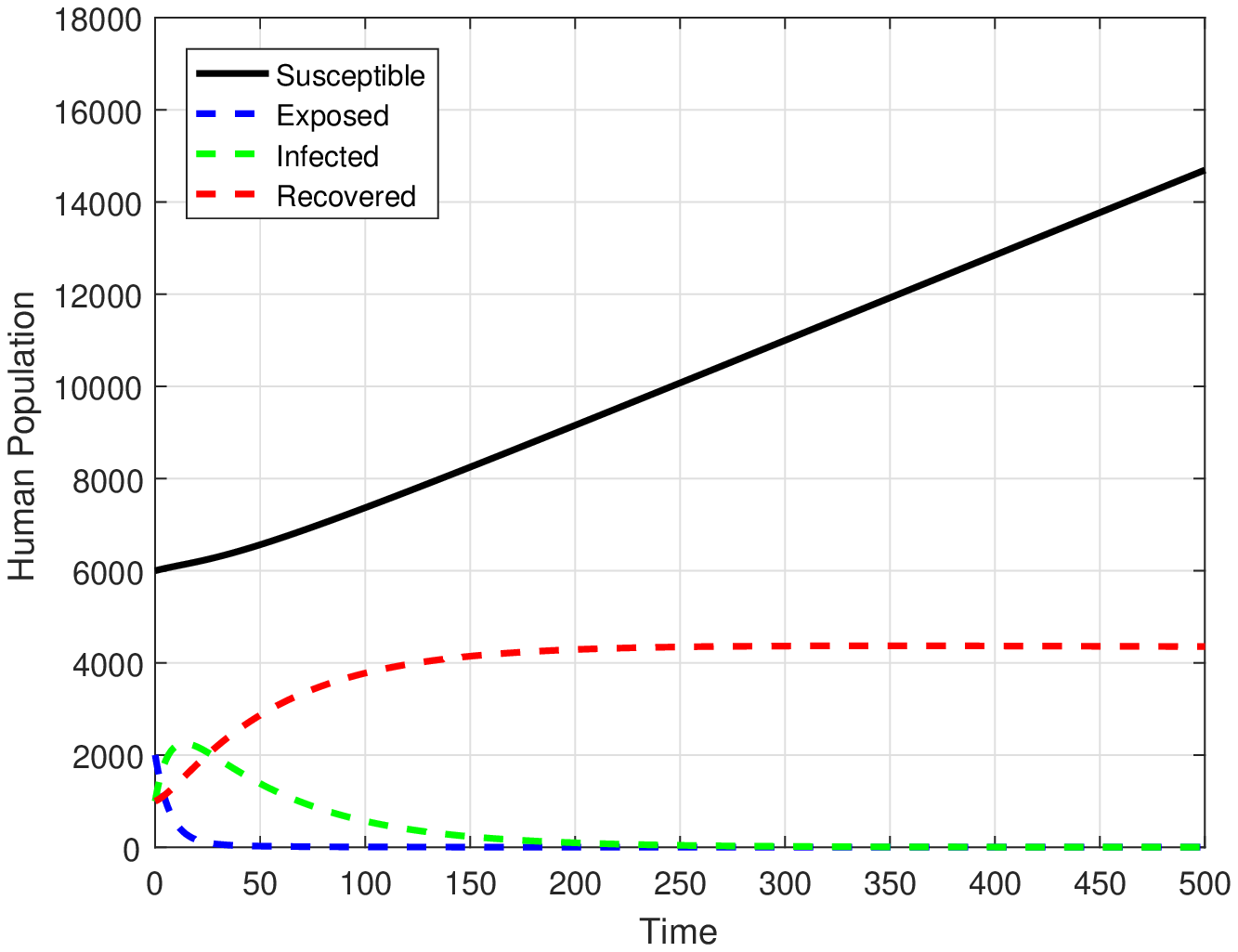}\label{b}}
\hspace{0mm}
\subfloat[$G_0=0, c_{vh}=1, c_{hv}=1, a_v=3.0425$]{\includegraphics[scale=0.5]{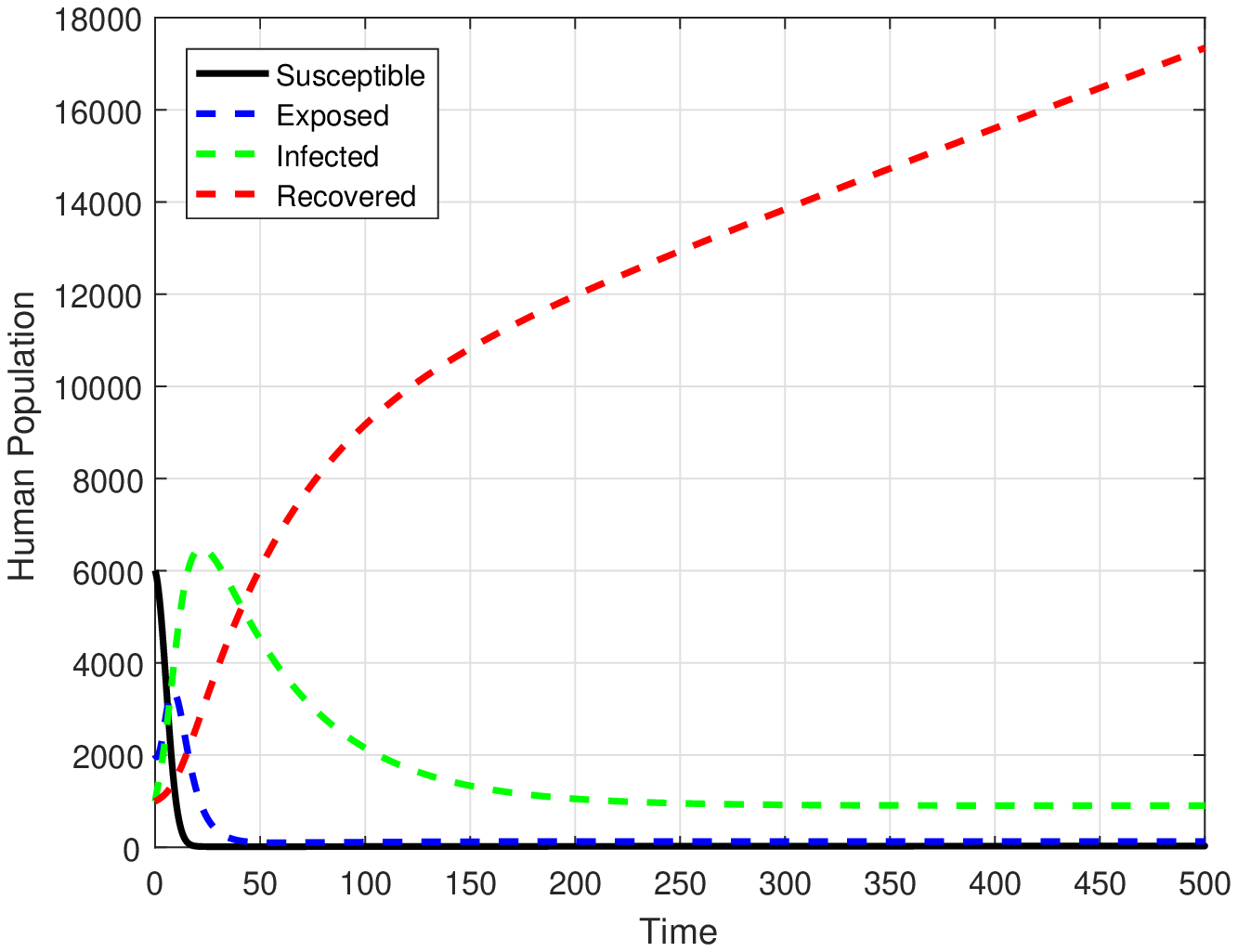} \label{c}}
\subfloat[$G_0=20,c_{vh}=1, c_{hv}=1, a_v=3.0425$]
{\includegraphics[scale=0.5]{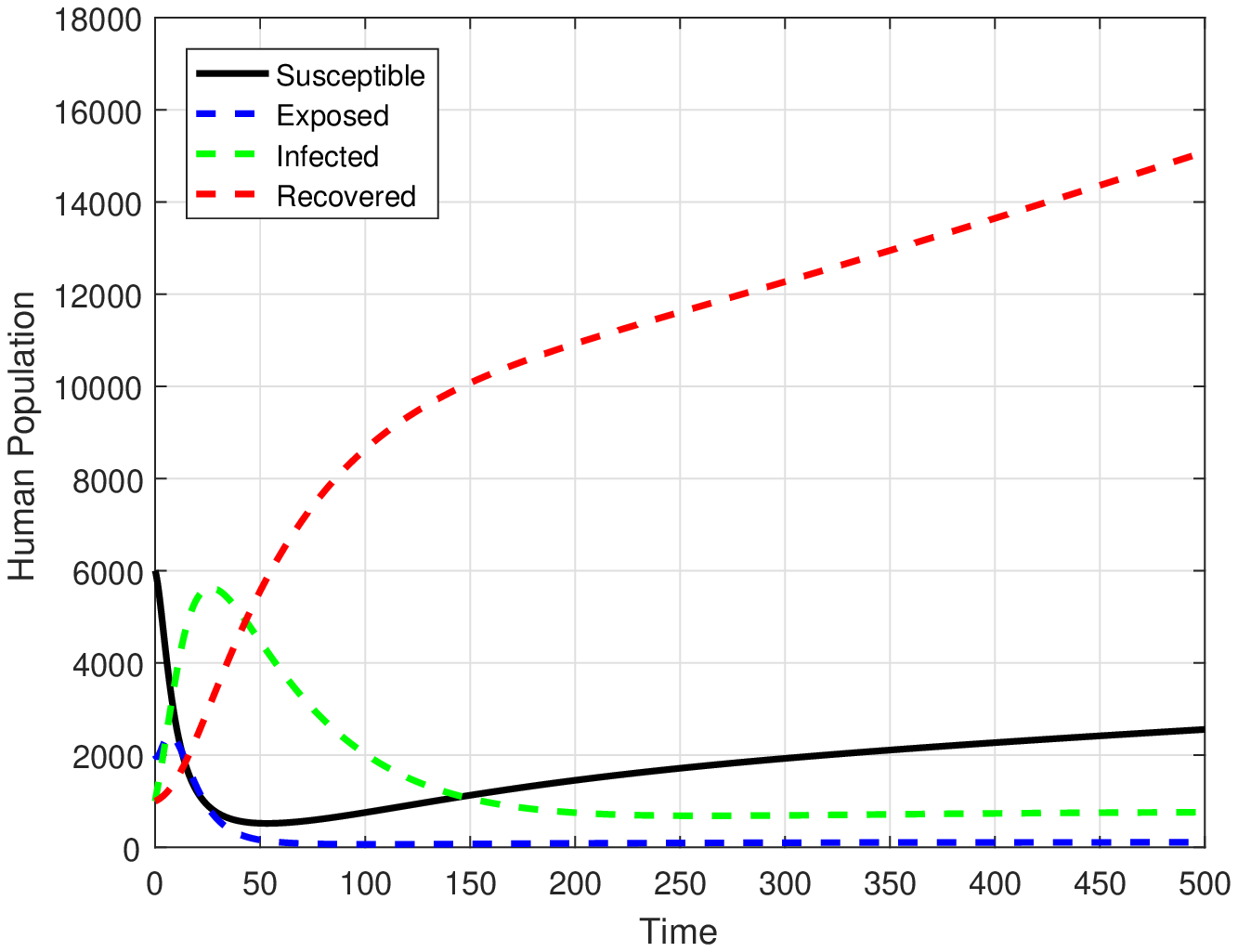}\label{d}}
\end{center}
\caption{Effect of presence of predator on the human population within 500 days}\label{plot1}
\end{figure}

\begin{figure}[h]
\begin{center}
\subfloat[$G_0=0, c_{vh}=0.2, c_{hv}=0.25, a_v=0.25$]{\includegraphics[scale=0.5]{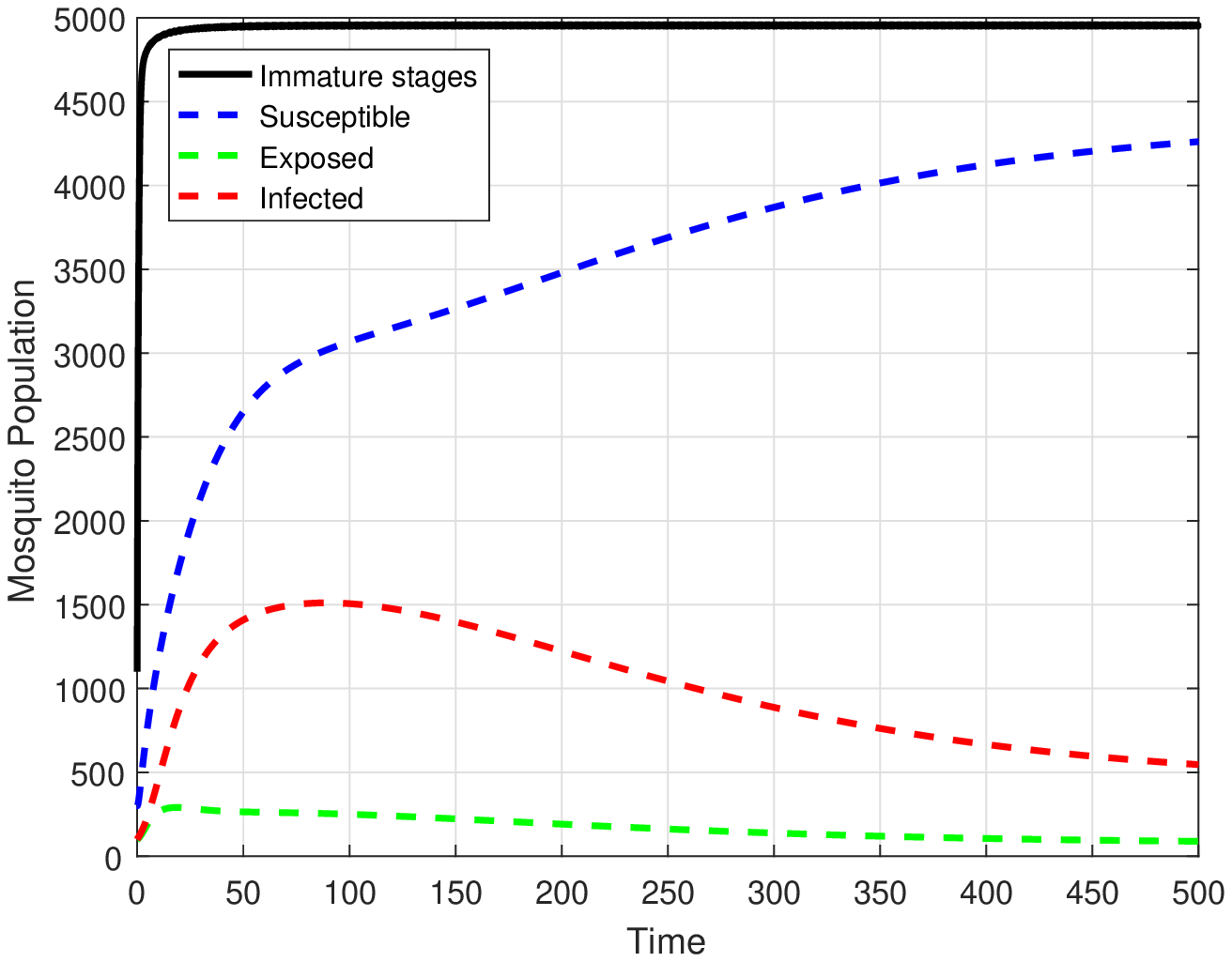} \label{e}}
\subfloat[$G_0=20, c_{vh}=0.2, c_{hv}=0.25, a_v=0.25$]{\includegraphics[scale=0.5]{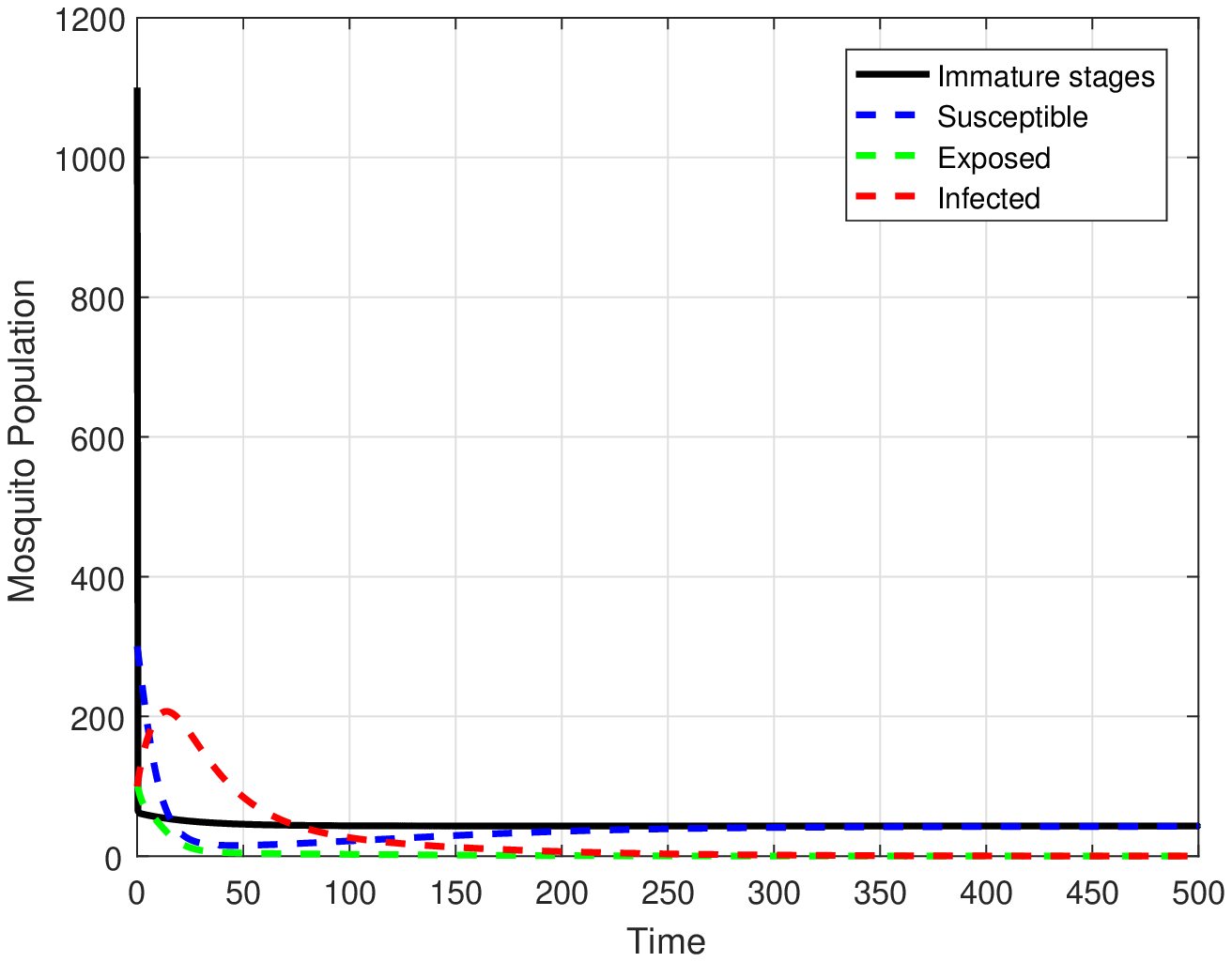}\label{f}}
\hspace{0mm}
\subfloat[$G_0=0, c_{vh}=1, c_{hv}=1, a_v=3.0425$]{\includegraphics[scale=0.5]{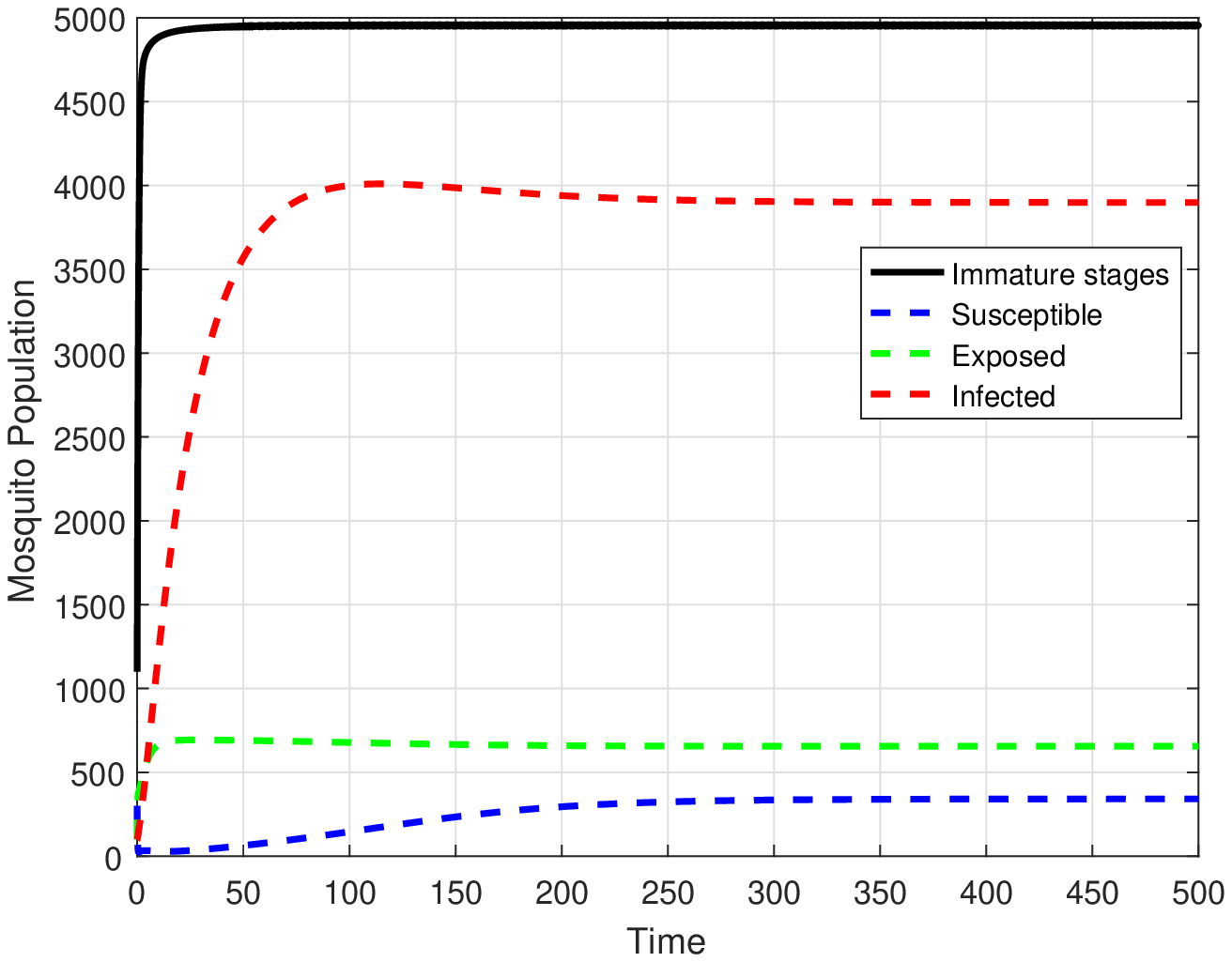} \label{g}}
\subfloat[$G_0=20, c_{vh}=1, c_{hv}=1, a_v=3.0425$]
{\includegraphics[scale=0.5]{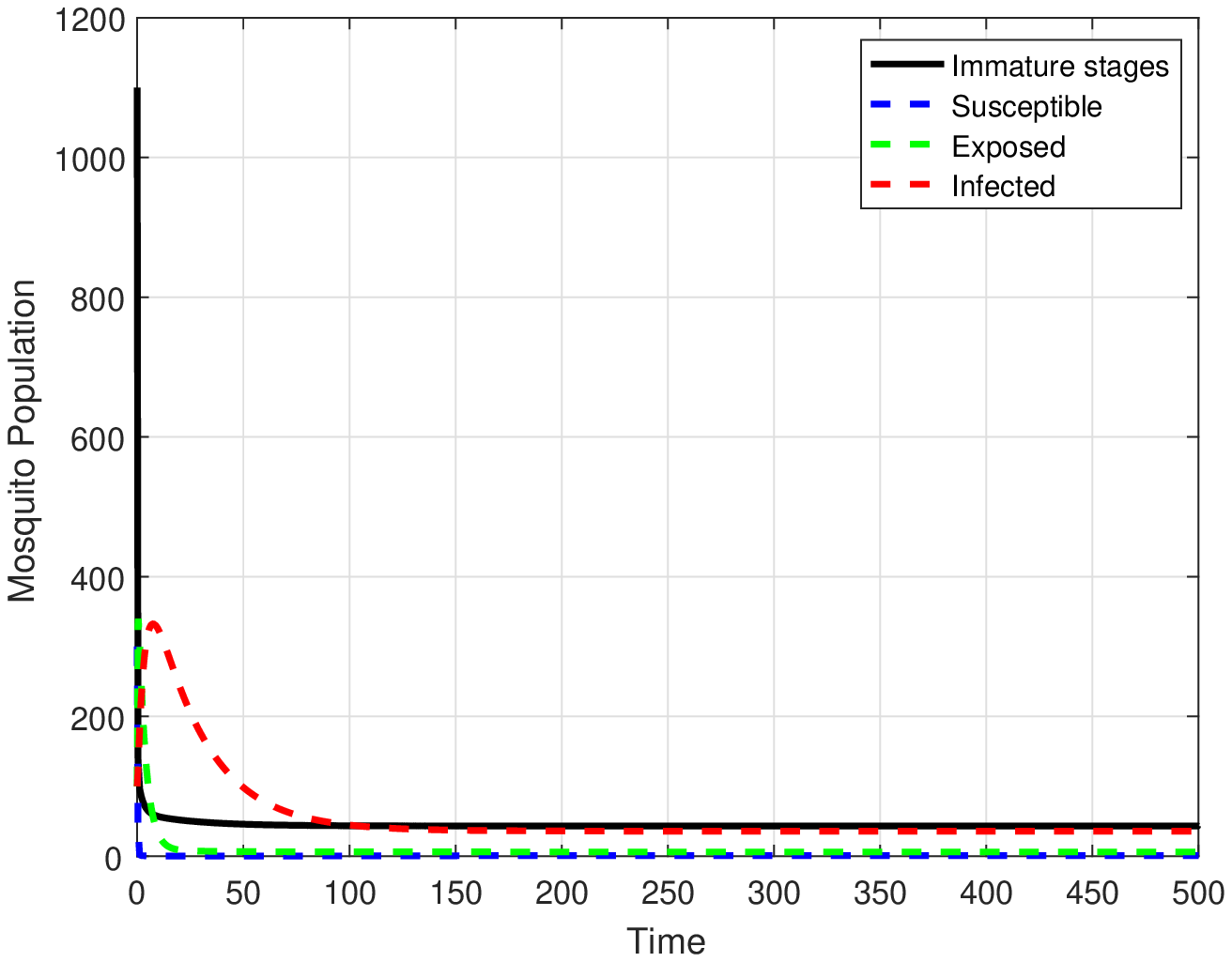}\label{h}}
\end{center}
\caption{Effect of presence of predator on the mosquito population within 500 days}\label{plot2}
\end{figure}

\begin{figure}[h!]
\begin{center}
\includegraphics[scale=0.75]{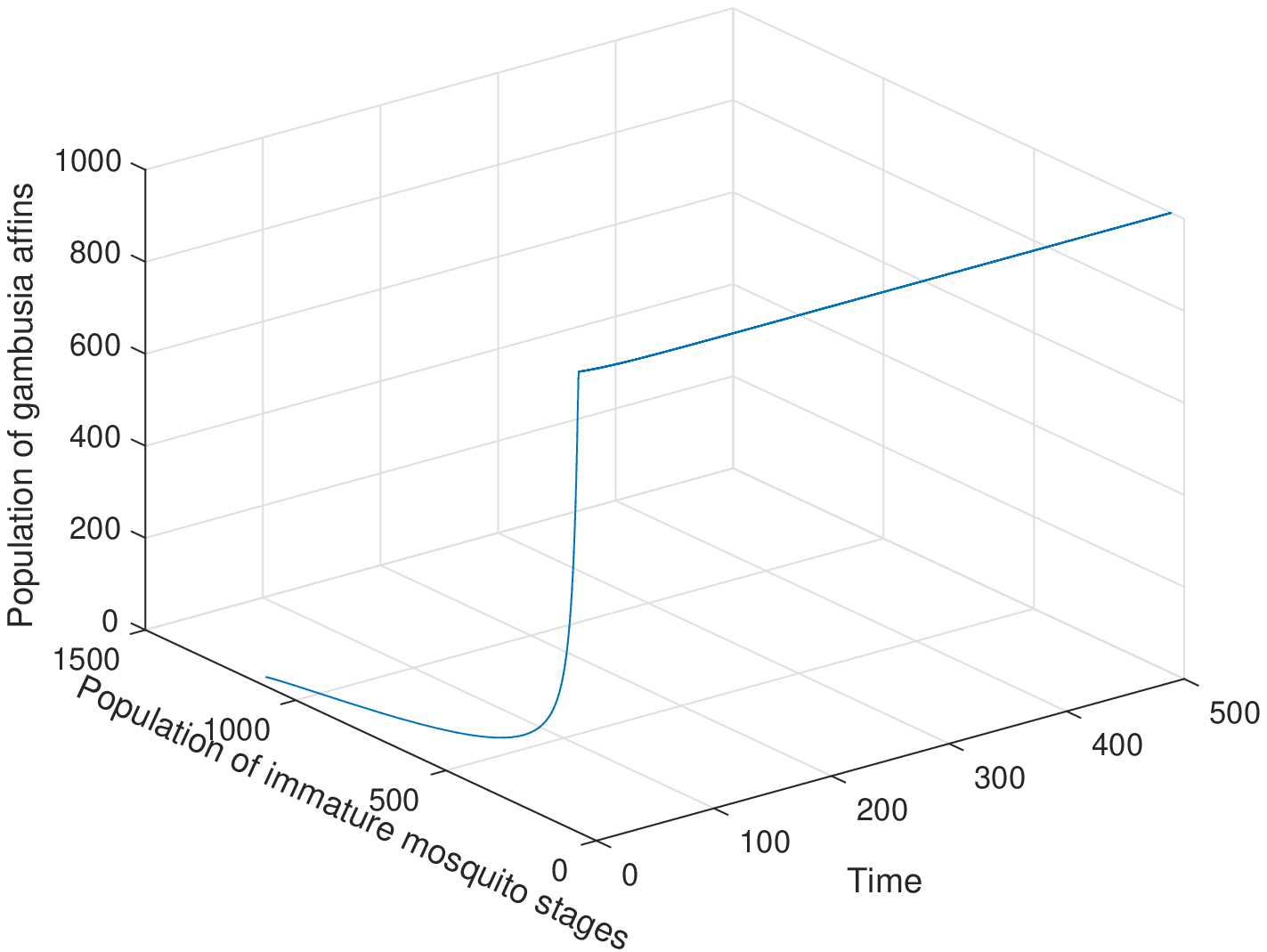}
\end{center}
\caption{Effect of introducing Gambusia Affinis in the system on immature stages of vector with time when $c_{vh}=0.2, c_{hv}=0.25, a_v=0.25$.}\label{plot3}
\end{figure}

This shows that practically pushing the reproduction number below unity is nearly impossible. This also explains the reason why controlling zika is such a huge problem. Among the factors in Table \ref{table:1}, we can most easily change the mosquito biting rate by being precautious like wearing full coverage clothes, using mosquito repellents, etc. If we consider $a_v$ to be a variable we can change, even then the value for basic reproduction number comes out to be $\mathfrak{R}_0=32.3028 a_v \sqrt{c_{vh}c_{hv}}$. Therefore, we conclude that relying on reproduction number solely is not of very much importance here. We try to see the effect of introducing the predator on the system in controlling the infection. The presence of predator would decrease the mosquito population which in turn would decrease the contact rates between host and vector.

For initial conditions we consider $S_{h_0}=6000, E_{h_0}=2000, I_{h_0}=1000, R_{h_0}=1000, m_{q_0}=1100, S_{v_0}=300, E_{v_0}=100, I_{v_0}=100$ and $G_0=0$ when predator is absent and $G_0=20$ when predator is present, that is, we introduce only 20 fishes initially.
For different reproduction numbers we assume $a_v=0.25, c_{vh}=0.2, c_{hv}=0.25$ and $a_v=3.0425, c_{vh}=1, c_{hv}=1$ which result in $\mathfrak{R}_0=1.81$ and $\mathfrak{R}_0=98.2814$ respectively. We can see clearly here that decreasing the reproduction number beyond unity is practically impossible. We plot graphs for the above two different values of reproduction number and try to analyze the effect of introducing the predator to the system.

Figure \ref{plot1} shows the effect of introducing the Gambusia Affinis on the different classes of the human population. We analyze that when $\mathfrak{R}_0=1.81$ and $G_0=0$, the infection persists even after the passage of 500 days while on introducing Gambusia Affinis infection disappears completely within 200 days only. Similarly, when $\mathfrak{R}_0=98.2814$ and $G_0=0$, infection is comparatively higher than the former case, and introducing the predator helps in decreasing the infection and increasing the susceptible population. We can visualize that when the basic reproduction number is very high, the infection in the population persists even on putting the predator into the system but it helps in decreasing the infection as well as increasing the susceptible population.

Figure \ref{plot2} depicts the effect of putting in predatory fish Gambusia Affinis into the system on the mosquito population. We can clearly see the inspiring effect of Gambusia on the mosquito population. When $\mathfrak{R}_0=1.81$, Gambusia helps in disappearing the infection completely as well as keeping the overall mosquito population in check including immature stages. When $\mathfrak{R}_0=98.2814$, again Gambusia decreases the infection to a very high extent and decreases the mosquito population in all the classes to a very high magnitude.  

Figure \ref{plot3} represents the relative change in the population of Gambusia Affinis and the population of immature stages with respect to time. We can see that by introducing merely 20 fishes into the system that can hold up to a maximum of 1000 fishes, the population of immature stages decreases drastically and the predatory fish reaches its maximum carrying capacity before the passage of 100 days. Hence, Gambusia Affinis can be a great tool to control the zika epidemic anywhere.

\section{Conclusion}\label{sec4}
 The system \eqref{eq:main} is a nine compartmental zika model where we analyze the effect of introducing the predatory fish Gambusia Affinis in decreasing the infection in the human population even when the reproduction number is not less than unity. We deduce that in the case of zika infection reducing the basic reproduction number below unity is practically extremely difficult. Therefore, we try to push the concerned mosquito threshold parameter below unity. Hence, in such a scenario Gambusia Affinis can be of great aid since till now Zika has no treatment or preventive vaccine. Moreover, we introduce Gambusia in very small man-made ponds that can hold up to a maximum of 1000 units of fish. We find that when the reproduction number is closer to unity but not less than unity, Gambusia Affinis can completely wipe off the infection from human as well as vector population in a maximum of 200 days period. Hence, together with the attempts to reduce mosquito biting rate and mutual contact rates between host and vector, introducing Gambusia Affinis can be of great help to eliminate the zika infection. We find that when the basic reproduction number is very high even then this predatory fish helps in increasing the susceptible population as well as keeping a check on the infected population. We also find that when there is no sexual transmission of zika then in absence of mosquito population or in absence of infection in the mosquito population, human zika infection can not persist. As it is already known that there have been reported a very few cases of sexual transmission of zika, if one can be precautious enough to not spread it sexually then we can get rid of zika infection by only getting rid of the mosquito population. Gambusia Affinis can be prominent support for this. Moreover, Gambusia Affinis is a fish breed that does not need much human attention and resources for growth and development, and using it in small ponds can help save the ecosystem balance too. Hence, the predatory fish Gambusia Affinis is of great importance in eradicating zika infection.
 
 Like all mathematical models, our model also has some limitations. We do not consider the sexual transmission of Zika infection which may be an important aspect and leave it for future study. Moreover, in the future predator can be introduced in larger ponds with harvesting to keep the population in control and inhibit it to cause ecological instability.



%
\section*{Conflict of interest}
The authors declare that they have no conflict of interest.


\begin{thebibliography}{99}

\bibitem{le} Global health observatory (gho) data. http://www.who.int/gho/mortality\_burden\_disease/life\_tables/situation\_ trends/en/ (2016). 

\bibitem{wb} Birth rate, crude (per 1,000 people). https://data.worldbank.org/indicator/SP.DYN.CBRT.IN (2017). 

\bibitem{incubation} World health organisation. https://www.who.int/news-room/fact-sheets/detail/zika-virus (2018). 

\bibitem{gambusiapredatoryuse} Gambusia aﬃnis (baird and girard, 1853). https://nas.er.usgs.gov/queries/factsheet.aspx?speciesID=846 (2019). 

\bibitem{aedesls} Life cycle of aedes aegypti. http://www.denguevirusnet.com/life-cycle-of-aedes-aegypti.html (2019). 

\bibitem{arias} Arias, J.H., Martínez, H.J., Sepúlveda-Salcedo, L.S., Vasilieva, O.: Predator-prey model for analysis of aedes aegypti population dynamics in Cali, Colombia. Int. J. Pure Appl. Math. 105(4), 561--597 (2015). 

\bibitem{gam2} Barnickol, P.G.: Food habits of gambusia aﬃnis from reelfoot lake, tennessee, with special reference to malaria control. Rep. Reelfoot Lake Biol. Sta. 5(4), 5–-13 (1941).

\bibitem{zika5} Binder, M., Pilyugin, S.S.: Stability analysis of a deterministic model of Zika/Dengue co-circulation. Int. J. Biomath. 12(4), 1950045, 39 (2019). 

\bibitem{zika4} Cai, Y., Ding, Z., Yang, B., Peng, Z., Wang, W.: Transmission dynamics of Zika virus with spatial structure-a case study in Rio de Janeiro, Brazil. Phys. A 514, 729–-740 (2019).


\bibitem{zika1} Cai, Y., Wang, K., Wang, W.: Global transmission dynamics of a Zika virus model. Appl. Math. Lett. 92, 190--195 (2019). 

\bibitem{gam7} Chandra, G., Bhattacharjee, I., Chatterjee, S., Ghosh, A.: Mosquito control by larvivorous ﬁsh. Indian J Med Res. 127(1), 13--27 (2008). 

\bibitem{chitnis} Chitnis, N., Hyman, J., Cushing, J. Determining important parameters in the spread of malaria through the sensitivity analysis of a mathematical model. Bull. Math. Biol. 70, 5 (2008), 1272--1296.

\bibitem{2overview} Dick, G.W.A.: Zika virus (ii). pathogenicity and physical properties. Trans R Soc Trop Med Hyg. 46(5), 521--534 (1952). 

\bibitem{Dkmbook} Diekmann, O., Heesterbeek, J.A.P.: {\it Mathematical epidemiology of infectious diseases}. Wiley Series in Mathematical and Computational Biology. John Wiley \& Sons, Ltd., Chichester (2000).

\bibitem{Dkm1} Diekmann, O., Heesterbeek, J.A.P., Metz, J.A.J.: On the deﬁnition and the computation of the basic reproduction ratio R0 in models for infectious diseases in heterogeneous populations. J. Math. Biol. 28(4), 365--382 (1990).

\bibitem{Dkm2} Diekmann, O., Heesterbeek, J.A.P., Roberts, M.G.: The construction of next generation matrices for compartmental epidemic models. J. R. Soc. Interface 7, 873--885 (2009).

\bibitem{geneticex} Evans, B., Kotsakiozi, P., Costa-da Silva A.L., e.a.: Transgenic aedes aegypti mosquitoes transfer genes into a natural population. Sci Rep. 9(1) (2019). 

\bibitem{gam3} Gerberich, J.B.: An annotated bibliography of papers relating to the control of mosquitoes by the use of ﬁsh. Am. Mid. Nat. 36(1), 87--131 (1946). 

\bibitem{gam4} Holland, E.A.: An experiment in the control of malaria in new ireland by distribution of gambusia aﬃnis. Trans. Roy. Soc. Trop. Med. and Hyg. 26(6), 529--538 (1933). 

\bibitem{di} Imran, M., Usman M.and Dur-e-Ahmad, M., et al.: Transmission dynamics of zika fever: A seir based model. Diﬀer Equ Dyn Syst 26, 1--24 (2017). 

\bibitem{ppm} Jana, S., Kar, T.: A mathematical study of a prey-predator model in relevance to pest control. Nonlinear Dyn 74(3), 667--683 (2013). 

\bibitem{gam6} Jordan, D.S.: Malaria and the mosquito ﬁsh. Scient. Amer. 134(5), 296--297 (1926). 


\bibitem{gam5} Jordan, D.S.: The mosquito ﬁsh (gambusia) and its relation to malaria. Smithsonian Rept. pp. 361–368 (1926). 


\bibitem{3overview} K.C., S.: Neutralizing antibodies against certain recently isolated viruses in the sera of human beings residing in eastafrica. J Immunol. 69(2), 223--234 (1952). 

\bibitem{gam1} Krumholz, L.A.: Reproduction in the western mosquitoﬁsh, gambusia aﬃnis aﬃnis (baird and girard), and its use in mosquito control. Ecological Monographs 18(1), 1--43 (1948). 

\bibitem{zika3} Nda\"{\i}rou, F., Area, I., Nieto, J.J., Silva, C.J., Torres, D.F.M.: Mathematical modeling of {Z}ika disease in pregnant women and newborns with microcephaly in {B}razil. Math. Methods Appl. Sci. 41(18), 8929--8941 (2018).

\bibitem{recovery} Paz-Bailey, G., Rosenberg, E., Doyle, K., Munoz-Jordan, J., Santiago, G., Klein, L., Perez-Padilla, J., Medina, F., Waterman, S., Gubern, C., et al.: Persistence of zika virus in body fluids–final report. N. Engl. J. Med. 379(13), 1234--1243 (2018). 

\bibitem{Pyke} Pyke, G.H.: A review of the biology of gambusia aﬃnis and g. holbrooki. Sci Rep. 15(4), 339--365 (2005). 


\bibitem{gam9} Roux, O., Robert, V.: Larval predation in malaria vectors and its potential implication in malaria transmission: an overlooked ecosystem service? Parasit Vectors 12(1) (2019). 

\bibitem{gam8} Shaalan, E., Canyon, D.: Aquatic insect predators and mosquito control. Trop Biomed. 26(3), 223--261 (2009). 

\bibitem{gam10} Shukla, A.S.: Biological control of vector of communicable disease in sidhi (mp) using larvivorous ﬁsh, gambusia aﬃnis. J. Math. Biol. 2(5) (2015). 

\bibitem{zika2} Srivastav, A.K., Yang, J., Luo, X., Ghosh, M.: Spread of Zika virus disease on complex network—a mathematical study. Math. Comput. Simul. 157, 15--38 (2019). 

\bibitem{zika7}  Wang, L., Zhao, H.: Dynamics analysis of a zika-dengue co-infection model with dengue vaccine and antibody-dependent enhancement. Physica A. 522, 248--273 (2019).

\bibitem{thpr} Yang, H., Macoris, M., Galvani, K., Andrighetti, M., Wanderley, D.: Assessing the effects of temperature on the population of aedes aegypti, the vector of dengue. Epidemiol Infect. 137(8), 1188--1202 (2009). 

\bibitem{biting} Yasuno, M., Tonn, R.: A study of biting habits of aedes aegypti in bangkok, thailand. Bull World Health Organ. 43(2), 319--325 (1970).

\end{thebibliography}


\end{document}